\documentclass{article}
\usepackage[utf8]{inputenc}
\usepackage{natbib}
\usepackage{amsmath}
\usepackage{mathtools}
\usepackage{amssymb}
\usepackage{hyperref}
\usepackage{cleveref}
\usepackage{caption}
\usepackage{subcaption}
\usepackage{soul}
\usepackage{fullpage}
\usepackage{afterpage}
\usepackage{graphicx}
\usepackage{fancyvrb}
\usepackage{booktabs}
\usepackage{placeins}
\usepackage{mymacros_10Sep19}
\usepackage{multirow}
\usepackage[version=4]{mhchem}
\usepackage[affil-it]{authblk}

\newcommand{\mreject}{r^\phi}
\newcommand{\mrecov}{r}
\newcommand{\flowr}{\dot{Q}}
\newcommand{\perm}{k}
\newcommand{\poros}{\epsilon}
\newcommand{\reac}{\xi}
\newcommand{\kin}{\mathcal{K}}

\newcommand{\diff}{D}

\newcommand{\openfoam}{OpenFOAM\textsuperscript{\textregistered}}

\renewcommand{\todo}[1]{\textcolor{red}{\textsc{#1}}}

\title{Mathematical modelling and numerical simulation of reverse-osmosis desalination}
\author[1]{Nicodemo Di Pasquale  \thanks{Corresponding author: nicodemo.dipasquale@brunel.ac.uk}}
\author[2]{Mayo Akele}
\author[3]{Federico Municchi}
\author[2]{John King} 
\author[2]{Matteo Icardi} 

\affil[1]{Department of Chemical Engineering, Brunel University London, Uxbridge, UB8 3PH, United Kingdom}
\affil[2]{Department of Mathematics, University of Nottingham, Nottingham, NG7 2RD, United Kingdom}
\affil[3]{Department of Mechanical Engineering, Colorado School of Mines, Golden, CO 80401, US}
\date{
\today\\
Special issue on Process Intensification}

\begin{document}

\maketitle

\section*{Abstract}

The reverse osmosis membrane module is an integral element of a desalination system as it determines the overall performance of the desalination plant. The fraction of clean water that can be recovered via this process is often limited by salt precipitation which plays a critical role in its sustainability.
In this work, we present a model to study the complex interplay between flow, transport and precipitation processes in reverse osmosis membranes, which together influence recovery and in turn process sustainability. A reactive porous interface model describes the membrane with a dynamic evolving porosity and permeability to capture the scaling and clogging of the membrane. An open-source finite-volume numerical solver is implemented within the \openfoam library and numerical tests are presented here showing the effect of the various parameters of the model and the robustness of the model to describe a wide range of operating conditions.




\section{Introduction}

The demand of freshwater has steadily increased over the last forty decades at a global level, mainly because of the increase in population and improving living standards which are leading to an expansion of irrigated agriculture and its human consumption.
In turn, the increase in consumption is straining the freshwater sources in their ability to supply the growing demand of water worldwide with almost two third of the total world population experiencing severe water scarcity during at least a part of the year \citep{Mekonnen2016,Jones2019}. 
The current freshwater sources are already overexploited, and even if better management is still needed to reduce the  misuse of such resources (e.g., wastewater treatments or waste reduction) \citep{Najid2022}, these solutions cannot still be enough to meet the future demand of freshwater. The ongoing climate change is expected to reduce the availability of freshwater because of the receding of glaciers with a subsequent important reduction of the flow in important rivers such as the Mekong Yellow, or Gange \citep{Shannon2008}.  

Almost 98\% of the total liquid water on the Earth is not available for the direct use or consumption, as it form the total water present in seas and oceans.  However, this last fact also means that if the exceedingly high saline content in seawater can be reduced or removed, we will have access to the largest source of freshwater sources, with which the required amount of freshwater could be delivered without straining the natural occurring 
 resources \citep{Jones2019}.
Therefore, there is a strong push into developing more efficient technologies for the desalination of seawater. Among 
the currently available technologies are thermal desalination and membrane processes \citep{Fritzmann2007,Subramani2015}. 
In thermal desalination, seawater is brought to evaporation through multi-effect distillation or multi-stage flash distillation \citep{Al-hotmani2021}, and the resulting vapor is subsequently condensed. In membrane technologies, a semi-permeable membrane is employed to separate (or filtrate) the solution of salt and water.  

Reverse Osmosis (RO), is a widely employed membrane technology for treating seawater and wastewater with salinity up to 70 g/l \citep{Hickenbottom2014} which, due to its relative simplicity and widespread diffusion has 
been among the main topics of research in membrane filtration \citep{Wardeh2008,Luo2019}.
One of the main challenges in RO is \emph{concentration polarisation} \citep{Kim2005}, which is the presence, over the
membrane, of a solute-rich boundary layer. Concentration polarisation can lead to solute precipitation and fouling,
significantly reducing the local permeability of the membrane with adverse effects on the permeation flux. When
the solution contains inorganic salts (such as sodium chloride NaCl or calcium sulfate CaSO$_4$ ), the resulting
accumulation of crystals on the membrane is also called \textit{scaling}. There is extensive experimental evidence indicating that scaling reduces the membrane performance over time \citep{Hu2014} and therefore,  scaling 
phenomena should play a major role in the mathematical modelling of RO systems.

Computational Fluid Dynamics (CFD) represents a powerful to analyse concentration polarisation and scaling, with  the earliest attempts dating back to more than two decades ago \citep{Hansen1998}. The membrane is usually included as a boundary condition where the flux of the solute is assumed equal to zero (complete retention). Previous studies have
considered dependence on the solute concentration of properties such as the osmotic pressure \citep{Hansen1998}, viscosity, density and diffusion coefficient \citep{Geraldes2001, Wiley2002} to include effects of
concentration polarisation on the membrane\citep{Wiley2002,Johnston2022}.
In particular, CFD simulations for membranes have been employed to analise different geometrical configurations, such as the spacer-filled channels. These include solid elements in the feeding flow to increase the shear stress on the surface of the membrane, which in turn increases the local mixing and mass transfer across the membrane \citep{Shakaib2007, Fletcher2004, Fimbres2010, Ghidossi2006,Lau2009,Santos2007,Koutsou2009,Ranade2006}.

However, the effects of scaling have not been directly included in CFD simulations.
In this work, we propose a mathematical model and a CFD solver for the analysis of the performances of a RO membrane in which we included the possibility for the solute in the feed flow to react (precipitate) on the membrane and therefore affecting the membrane permeability.
The paper is organised as follows. We describe the mathematical framework for the analysis of the membrane, highlighting how the chemical reactions can be accounted for in the model. We then discuss the implementation of our model into the widely used open source finite volume library \openfoam  and we show some results for some typical situations, showing the applicability of the whole framework. We then draw some conclusions and we outline the possible extension of the model. 

\section{Model}\label{sec:model}

\begin{figure*}
   \includegraphics[width=\textwidth]{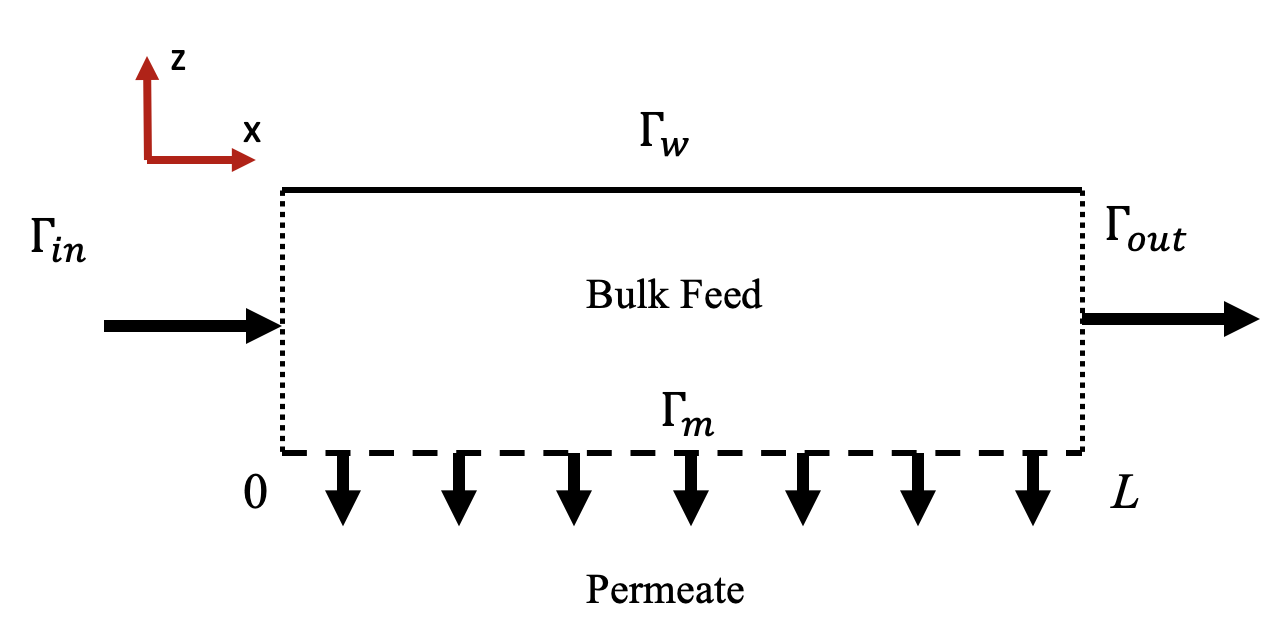}
 \caption{2-dimensional domain considered in this work
 }
    \label{fig:geometry}
\end{figure*}

In this work, we approximate a rectangular 3D membrane module as a 2D channel illustrated in \cref{fig:geometry}. This is a common practice employed in other recent CFD studies 
\citep{Johnston2022}. While the configuration we considered can be easily extended to more complex geometries, the focus here is to develop a complete mathematical model for the polarization and scaling of the membrane by giving a proof of concept for a general-purpose numerical solver. Our main goal is to show the mechanisms governing the evolution of the solute at the membrane interface, and a 2D channel geometry allows us to focus on this task.

Let us therefore assume a rectangular domain $\Omega \equiv (0,L) \times (0,H)$ with boundary $\Gamma = \partial \Omega$ subdivided in three different regions:
\begin{align}
    \Gamma_m & = (0, L) \times \{0\} \quad \text{(membrane)} \nonumber  \\ 
    \Gamma_{in} & =   \{0\} \times (0, H) \quad \text{(inlet)} \nonumber \\  
    \Gamma_{out} & =   \{L\} \times (0, H) \quad \text{(outlet)}
\end{align}
where $\Gamma_m$ represents the membrane boundary, $\Gamma_{in}$ is the inlet boundary and $\Gamma_{out}$ is the outlet boundary. With $\Gamma_w$ we represent the remaining part of the boundary constituted by solid boundaries so that $\Gamma_w = \Gamma \setminus (\Gamma_{in} \cup \Gamma_{out} \cup \Gamma_m)$, as shown in \cref{fig:geometry}.
In this study, the permeate flow is not explicitly modelled, therefore the membrane represents a boundary condition for the problem.

\begin{figure*}[!t]
    \centering
   \includegraphics[width=12cm]{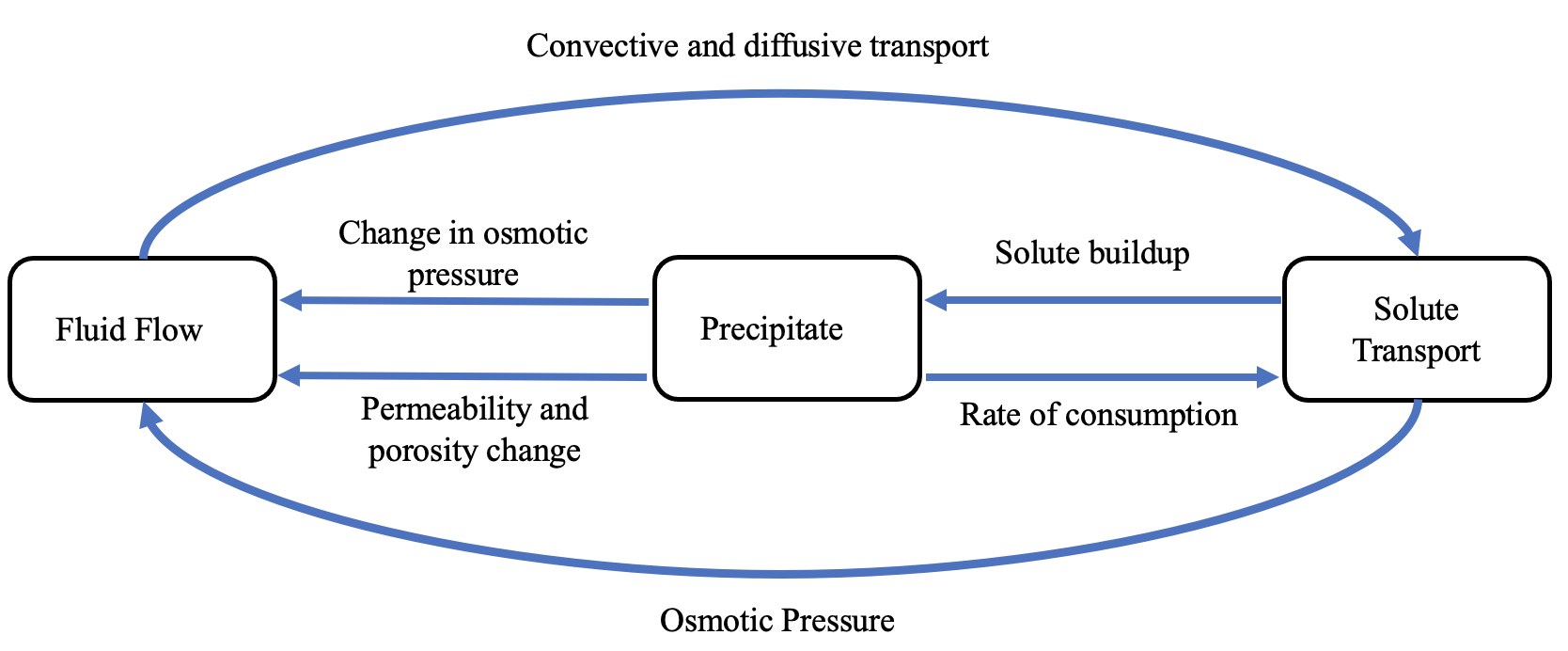}
 \caption{Sketch of the different mechanisms and models considered in this work.  }
    \label{fig:coupling}
\end{figure*}

The analysis of the filtration process requires the simultaneous solution of the flow field coupled with the transport of dissolved chemical species, which can involve one or more chemical reactions. Such chemical
reactions can lead to solute precipitation, modifying the permeability and porosity of the membrane, thus causing a variation in the osmotic pressure and the final yield of the permeate. We summarised the different mechanism and their interdependence in \cref{fig:coupling}.

Our mathematical model is composed of three different and intertwined parts that need to be simultaneously considered to obtain the overall behaviour of the membrane which we are now going to describe in detail.

\subsection{Flow modelling}

The mixture is described as an incompressible Newtonian fluid described by the Navier-Stokes equations:
\begin{subequations} 
\label{eq:NS}
\begin{align}
     \nabla \cdot  \bfu & =0  & & \mbox{in}\:\:\: \Omega, t>0 \label{eq:cont} \\
    \rho \sqp{\pard{\bfu}{t} + \cip{\bfu\cdot\nabla\bfu}} & = -\nabla p + \mu \nabla^2\bfu & & \mbox{in}\:\:\: \Omega, t>0 
\end{align}
\end{subequations}
where $\rho$ is the density of the fluid $\bfu$ is the velocity vector with components $(u, v)^T$ , $p$ is the fluid pressure and $\mu$ is the dynamic viscosity, $\nabla$ is the gradient operator. We impose the following boundary and initial conditions on the system shown in \cref{eq:NS} for the velocity and pressure:
\begin{subequations}
\label{eq:allBCsNS}
\begin{align}
    \bfu & = 0 & & \mbox{in}\:\:\: \Omega, t=0  \label{eq:pc1}\\
     \bfu & = 0,\:\:\   p=P_{in}  & & \mbox{on}\:\:\: \Gamma_{in} \label{eq:pc2} \\ 
     \nabla u\cdot \bfn  & =0, \:\:\   p=P_{out} \:\:  & & \mbox{on}\:\:\: \Gamma_{out}  \label{eq:pc3} \\     
     \bfu & = 0\:\: & & \mbox{on}\:\:\: \Gamma_w  \label{eq:pc4} 
\end{align}
\end{subequations}
Finally, the membrane is modelled as a dynamic Dirichlet condition for the velocity $v$ orthogonal to the membrane is obtained from the Darcy law as: 
\begin{align}\label{eq:bcmembr}
    u = 0, \:\: v=-\frac{\perm(\Delta p - \Delta \pi)}{\ell \mu} & & \mbox{on}\:\:\: \Gamma_m
\end{align}
where $\mu$ is the fluid viscosity, $\ell$ is the membrane thickness, $k$ is the membrane permeability, $\Delta p = p_m - p_p$ is the pressure difference between the feed side of the membrane, $p_m$, calculated at the membrane boundary $\Gamma_m$ and the permeate side of the membrane, $p_p$, whereas $\Delta \pi$ is the osmotic pressure gradient. 
This last equation shows that the flux through the membrane is proportional to the difference between the applied and osmotic pressure differentials.

The pressure gradient $\Delta p=p_m - p_p$ is the difference between the pressure on the feed side of the membrane $p_m$ calculated at the membrane boundary $\Gamma_m$, and the permeate side of the membrane $p_p$. The osmotic pressure difference is defined as follows:
\begin{equation}\label{eq:osmp}
    \Delta \pi = \pi - \pi_p = \pi(\phi_1,\ldots,\phi_N) - \pi(\phi^p_1,\ldots,\phi^p_N)
\end{equation}
where the feed osmotic pressure $\pi$ is a function of the $N$ ions concentrations, $\phi_i$, $i=1,\ldots,n$, whereas the permeate osmotic pressure $\phi_p$ is function of the $N$ ions concentrations in the permeate,  $\phi^p_i$, $i=1,\ldots,n$.

As common in membrane applications \citep{Linares2014} the osmotic pressure is computed using the Van't Hoff model \citep{VantHoff1888}:
\begin{equation}\label{eq:pitzer}
    \pi(\phi_1,\dots,\phi_n) = RT\varphi \sum_{i}^{N}\phi_i
\end{equation}
where $\varphi$ is the osmotic coefficient, $R$ is the gas constant and $T$ is the temperature we can rewrite \cref{eq:osmp} as:
\begin{equation}\label{eq:deltapi}
     \Delta \pi = RT \varphi \sum_{i}^{N}\cip{\phi-\phi^p_i} = RT  \varphi \sum_{i}^{N}\mreject_i\phi_i
\end{equation}
where we used the definition of membrane rejection of the species $i$, 
\begin{equation} \label{eq:rej}
     \mreject_i = 1- \frac{\phi^p_i}{\phi_i},\:\: i=1,\ldots,N \, .
\end{equation}
The Van't Hoff equation assumes a linear relation between concentration and osmotic pressure and is more accurate for low concentration of solutes. Different formulations were proposed to take into account the non-linear behaviour of solutions, which consider the activity of the solvent \citep{Khraisheh2019}, calculated using the Pitzer equation for the electrolyte solutions \citep{Pitzer1973}. As a proof of concept, we will consider here the simplified model since the most complex behaviour can be straightforwardly added to the model and the CFD code we will present in the next section.
The membrane rejection of the species $i$ expresses the amount of solute rejected by the membrane (and therefore not present in the permeate) as a fraction of the initial quantity. In this work, we assume $\mreject_i=1$ for every ion species, which corresponds to complete rejections of the ions at the membrane. In this case, the concentration on the permeate side and the permeate osmotic pressure, $\pi_p$, are both equal to $0$ and therefore:
\begin{equation}
     \Delta \pi =  RT  \varphi \sum_{i}^{N}\phi_i\, .
\end{equation}
In the literature, the two most popular models proposed to describe the solute-solvent solute transport through the membrane are the solution-diffusion model \citep{Merten1963, Lonsdale1965,Wijmans1995} and the Spiegler-Kedem model \citeauthor{Spiegler1966}. The former expresses the flow through the membrane $\dot{J_v}$ as: 
\begin{align}\label{eq:diffmod}
    \dot{J_v}=A(\Delta p - \Delta \pi)
\end{align}
noting that in our notation  $v=\frac{\dot{J_v}}{A(\Gamma_m)}$ where $S_m$ is the membrane surface; while for the latter we have  
\begin{align}\label{eq:spiegmod}
    \dot{J_v}=\frac{1}{R_m\mu}(\Delta p - \sigma \Delta \pi)
\end{align}
with the same identification for $\dot{J_v}$ (i.e., $v=\frac{\dot{J_v}}{A(\Gamma_m)}$)   where $R_m$ is the membrane resistance and $\sigma$ is the reflection coefficients which measure the impermeability of the membrane to the solutes. $\sigma=1$ indicates a membrane completely impermeable to solutes and will be the one considered for this work.

Notice that in \cref{eq:bcmembr} we use the Darcy law to rewrite the so-called water permeability of the membrane $A$ which appear in the equation in terms of the permeability and thickness of the membrane and the viscosity of water as $A=\frac{S_m\perm}{\ell\mu}$. 
Using the same reasoning for \cref{eq:spiegmod} we obtain that  $R_m=\frac{\ell}{S_m\perm}$, that is to say, the membrane resistance is inversely proportional to the permeability. 
The identification of the Darcy-related terms with the water permeability $A$ or membrane resistance has two main advantages. Firstly, we can connect our analysis with the membranes available commercially which are described in terms of water permeability or membrane resistance. This last fact allows us to choose the range of the parameters we are considering which are appropriate for the description of real membranes. Secondly, using the Darcy derived version of the equation for the flow through the membrane allows us to include more detailed mechanisms in the model which can take into account more complex phenomena, such as chemical reactions and depositions of solids on the membrane as we will show in more details in the next sections.
One of the ways considered in the literature to include fouling and polarisation of the membrane is to define such contribution as additional resistance terms to be added to $R_m$ in \cref{eq:spiegmod} (see e.g., \citep{Silva2011,Lee1998membrane,Yeh2002}). These additional terms must be derived from experiments or empirical correlations. In contrast, our formulation leverages the well-established theory of porous media to include such effects directly in the determination of the permeability $k$.


The performance of the membrane can be evaluated by using the recovery $\mrecov$, defined as the ratio between the permeate flow rate, $\flowr_p$, and the feed flow rate, $\flowr_f$:
\begin{equation}\label{eq:recov}
    \mrecov=\frac{\flowr_p}{\flowr_f} = \frac{vA_p}{U_{in}A_{in}} = \frac{vL}{uH}
\end{equation}
where we used the definition of flow rate through a surface, $AU/\ell t$, to calculate the flux through the inlet (subscript $in$) and the membrane (subscript $m$) and then we replaced the area in \cref{eq:recov} with their geometrical expression related to the domain we are considering: $A_{in}=H\times Z$, $A_{p}=L\times Z$.

\subsection{Solute Transport}


In membrane processes, flow and solute transport are tightly coupled through the boundary condition on
the membrane given by equation \cref{eq:bcmembr}.  The transport equation for the concentration $\phi_i$ of the $i-$th ion species in the bulk liquid is given by:
\begin{subequations}
\label{eq:transp}
\begin{align} 
    \pard{\phi_i}{t}  & + \nabla\cdot(\bfu_i \phi_i) = \nabla \cdot \cip{\diff \nabla \phi_i} + \reac_{B}^i & & \mbox{in}\:\:\: \Omega, t>0 \label{eq:solute} \\
    \phi_i & = 0  & & \mbox{in}\:\:\: \Omega, t=0 \\
    \phi_i & = \phi_{in}  & & \mbox{on}\:\:\: \Gamma_{in} \\    
    \nabla \phi_i & = 0  & & \mbox{on}\:\:\: \Gamma_{out} \cup \Gamma_{w} \\
     \dot{J}_i & = v\phi_i - \diff_i \pard{\phi_{i}}{y} = \reac_{M}^i + \reac_{P}^i    & & \mbox{on}\:\:\:  \Gamma_m  \label{eq:BC}
\end{align}
\end{subequations}
where $\bfu_i$ is the velocity of the $i$-th chemical species, $D_i$ is its diffusion coefficient, $\reac$ is a rate term that represents different mechanisms of depletion of the solute, identified by the subscripts $B$ for the chemical reaction in the bulk,  $M$, for chemical reaction at the membrane, and $P$ for the fraction of the solute that crosses the membrane. The amount of the species $i$ at the membrane, can either precipitate on the membrane following a chemical reaction or go through the membrane in the permeate. The sum of these two mechanisms must equal the flux of the $i$-th chemical species at the membrane $\dot{J}_i$. We  can therefore assume that each of these mechanisms corresponds to a fraction of $\dot{J}$, and hence that $ \dot{J} = \kappa_M \reac_M + \kappa_P \reac_P$. The coefficients $\kappa$ have the property that:
\begin{equation}
    \kappa_M + \kappa_P = 1
\end{equation}
If only superficial reaction is present, then $\kappa_R=1$ and $\kappa_P = 0$, while in the case there is no superficial reaction and the solute crosses the membrane then $\kappa_P=1$ and $\kappa_R=0$. 
 
In this work we will make the following assumptions on the behaviour of the system: $1.$) The precipitation reaction is irreversible and only occurs on the membrane surface, (i.e., $\reac_B=0$) $2.$) The transport processes are similar for the all the salt ions and the diffusion coefficients are independent of concentration, $3.$) The salt ions have the same velocity as the fluid (i.e., $\bfu_i=\bfu$ for every $i$), $4.$) The porous medium is homogeneous, $5.$) there are only surface reactions at the membrane (i.e., $\kappa_R$ = 1).




\subsection{Chemical kinetics}\label{sec:chemkin}

Following the principles of mass action, the dynamic behaviour of chemical systems with $n$ components involved in $m$ reactions, can be described by a set of first order differential equations with time as the independent variable:
\begin{align}
	\totd{\phi_1}{t} & = f_1(\phi_1,\phi_2,\ldots,\phi_n,t) \nonumber \\
	\totd{\phi_2}{t} & = f_2(\phi_1,\phi_2,\ldots,\phi_n,t) \nonumber \\
 \vdots &\nonumber \\
	\totd{\phi_n}{t} & = f_n(\phi_1,\phi_2,\ldots,\phi_n,t) \nonumber \\
\end{align}
where $\phi_i(t)$, $i=1,\ldots,n$ denotes the volume molar concentration of chemical species $X_i$ at time $t$. The dynamics of the reaction network can be conveniently written in matrix form using the formalism developed in \citep{Chellaboina2009}:
\begin{equation}
	\totd{\phi_i}{t}  = \reac^i_R = (A - B)^T K \phi^A(t), \:\: \phi_i(0)=\phi_{0,i},\:\:t\geq0
 \end{equation}
where $K = \mbox{diag}(\kin_1, \ldots,\kin_m)$ is the diagonal matrix which contains as elements the reaction kinetics $\kin_j$, $j=1,\dots,m$ and $\phi_0$ is the initial concentration. $A$ and $B$ are the $m\times n$ matrices having in each entry the stochiometric coefficients of the reactants and products respectively and $\phi^A(t)$  is the matrix obtained by replacing each element of $A$ with $\phi_i^{a_{lp}}$ where $a_{kj}$ is the element of $A$ in the $l$-th row and $p$-th line.
The first term in the boundary conditions in \cref{eq:BC} takes the form
\begin{equation}
	\reac_R^i =  \at{\sqp{(A - B)^T K \phi^A(t)}}{i} \ell
\end{equation}
where the notation $\at{\sqp{\cdot}}{i}$ stands for the $i$-th component of its ($n\times 1$ vector) argument.

For a generic first order reaction $X_1+X_2 \rightarrow X_3$ with kinetic constant $K$, we have therefore
\begin{align}
	\totd{\phi_1}{t} & = -\kin \phi_1 \phi_2  \nonumber \\
	\totd{\phi_2}{t} & = -\kin \phi_1 \phi_2 \nonumber \\
	\totd{\phi_3}{t} & = \kin \phi_1 \phi_2	 \nonumber \\
\end{align}
and 
\begin{align}
	\reac_R^1 & = -\kin \phi_1 \phi_2 \ell  \nonumber \\
	\reac_R^2 & = -\kin \phi_1 \phi_2 \ell \nonumber \\
	\reac_R^3 & = \kin  \phi_1 \phi_2 \ell \nonumber \\
\end{align}
An example of a reaction of this type important in membrane operation is the formation of calcium carbonate, through the reaction \citep{Warsinger2015}:

\begin{equation}
\ce{Ca2+ + 2HCO3- -> CaCO3 + CO2 + H2O}    
\end{equation}

\noindent In fact, the scaling caused by the precipitation of calcium carbonate limits the operating condition of desalination systems for brackish, groundwater, and seawater \citep{Waly2009}.

One effect to be considered when a chemical reaction is involved is that the reaction between salt ions and subsequent precipitation of minerals, often alters the membrane properties, such as porosity and permeability. As crystals grow, it is expected that the permeability $\perm$, in equation (see \cref{eq:bcmembr}) and the porosity, $\poros$, will decrease, reducing the liquid flow through the membrane \citep{Steefel2005}. To account for the modification of the porosity and permeability as a result of mineral precipitation, we employ the Kozeny-Carman model \citep{Hommel2018}. This model allows us to  quantify the porosity-permeability relations and estimate the resulting changes.

According to the Kozeny-Carman model, the change in permeability is calculated by relating the current permeability, $\perm$, based on the current porosity $\poros$, to the initial permeability $\perm_0$ corresponding to the initial porosity $\poros_0$ \cite{Hommel2018}. These equations consequently follow the form below
\begin{equation}\label{eq:varporandper}
    \frac{\perm}{\perm_0}=\frac{f(\poros)}{f(\poros_0)}.
\end{equation}
Thus we can describe the permeability evolution using the following power law \citep{Hommel2018}
\begin{equation}\label{popem}
    \frac{\perm}{\perm_0} = \frac{(1-\poros_0)^2}{(1-\poros)^2} \cip{\frac {\poros}{\poros_0}}^3,
\end{equation}
where $\perm$ is the current permeability, $\poros$ is the current porosity, $\perm_0$ is the initial permeability and $\poros_0$ is the initial porosity. 
The rate at which porosity reduction occurs is given by \citep{huo2019,noiriel2004}:
\begin{equation}\label{eq:por}
    \epsilon = \epsilon_o - \frac{V_s}{\ell} \mathop \int \limits_{t_0}^{t} \sum_j \reac^j_R dt,
\end{equation}
where $\epsilon_o$ denotes the initial porosity at $t_0$, $V_s$ is the molar volume of solid precipitate in m$^3/$mol and $r(t)$ is the rate of precipitation in mol$\cdot \mbox{m}^{-2} \cdot \mbox{s}^{-1}$ and the index $j$ runs over all the possible reactions in the system. 
For the first order reaction we are considering here the \cref{eq:por} simplify as:
\begin{equation}\label{eq:prec}
    \epsilon = \epsilon_o - \frac{V_s}{\ell}\mathop \int \limits_{t_0}^{t} \left(\kin \phi_1\phi_2\right)\de t \,\,.
\end{equation}

\par We should again observe the inter-dependencies and feedback mechanisms between fluid flow, transport and reaction. Namely, in \cref{eq:prec} we see the precipitation reaction leads to a change in porosity which in turn affects the permeability in equation \cref{popem}. The change in permeability impacts the flow via the fluid velocity in \cref{eq:bcmembr}. This, in turn, alters the solute concentration distribution via \cref{eq:solute}, which ultimately impacts the rate of precipitation again via \cref{eq:BC}. Moreover, from \cref{eq:prec} we can observe that the variation of the porosity is proportional to the kinetic reaction constant. This last fact, in turn, simplifies the predictions for the clogging of the membrane based on the reactions in the systems. We can expect that if there are two chemical reactions in our systems, the first one 10 times slower than the second, then the clogging caused by the products of the second reaction will take 10 times the time needed by the products of the first reaction to produce the same amount of clogging.  
One of the strengths of our models is that allow these kinds of qualitative analyses even without actually solving the equations.

\section{Numerical discretisation}

The equations presented in \cref{sec:model} are solved using the open-source finite volume OpenFOAM 8.0 library and the code is available open-source \citep{membraneFoam}. In order to solve the equation of motion alongside the reaction at the membrane we developed a new solver for OpenFOAM called \textsf{binaryReactionFoam} which is based on two widely used solvers, \textsf{pimpleFoam} and \textsf{scalarTransportFoam}· The former is a transient solver for incompressible flows based on the PIMPLE algorithm while the latter is a concentration transport solver using a user-specified velocity field.  The solver also includes the possibility of modelling solid precipitation in the fluid and a multiphase flow model for the
solid particles. The most important element in the computational framework, however is represented by the new boundary conditions implemented to model the membrane. The \textsf{membraneVelocity} boundary conditions impose the fluid velocity based on the fluid pressure, the permeate pressure, and the membrane properties. These are updated in time by linking this boundary condition to the one for the scalar concentration, named \textsf{binaryReaction}, which solves for the solid precipitation at the boundary and therefore updates the membrane permeability. The equations and boundary conditions are coupled iteratively through Picard (fixed point) iteration (through the \textsf{PIMPLE iterations})
until convergence, making the whole model fully implicit.

We simulate a two-dimensional rectangular channel with height $h=0.003$ m and length of $L=0.02$ m discretised on a mesh composed by 600$\times$200  cells. 
The following discretisation schemes (we direct the reader to the OpenFOAM user guide \citep{Openfoam2019} for a detailed description of each scheme) are used to discretise the equations:
\begin{itemize}
    \item advective fluxes (\textsf{divSchemes Gauss vanLeer}) are computed at the faces and the variables interpolated with a Total Variation Diminishing scheme;
    \item gradient terms (\textsf{gradSchemes Gauss linear}) are approximated with central differencing;
    \item surface normal gradients for diffusive fluxes (\textsf{snGradSchemes orthogonal}) are approximated with central differencing (the grid is in fact orthogonal and does not need any correction to ensure second order accuracy);
    \item time derivatives (\textsf{ddTSchemes backward}) are approximated with third order implicit backward Euler scheme,
\end{itemize}
We specified a fully developed velocity profile at the inlet:
\begin{equation}
    u(0) =6 u_{av} \frac{y}{h}\cip{1-\frac{y}{h}} 
\end{equation}
where $u_{av}$ is the average velocity along the channel.
By specifying the velocity profile at the inlet we need only to specify the pressure at the outlet (the value of which is given in \cref{tab:consts}). The pressure of the permeate through the membrane is assumed constant along the length of the membrane and put equal to zero. For longer membranes, this assumption is no longer valid and the permeate flow needs to be modeled explicitly (with 1D or 2D models). This will be the subject of future extensions of our framework.

The initial value of the permeability we consider in the calculations is $k=10^{-16}$ m$^2$. Using the expression for the water permeability obtained from \cref{eq:bcmembr} and \cref{eq:diffmod}: $A=\frac{S_m\perm}{\ell\mu}$, and the viscosity of water at room temperature $\mu=10^{-3}$ Pa s, and the surface of the membrane, $S_m=2\cdot 10^{-6}$ m$^2$ for the channel configuration and $\ell=10^{-7}$m , we obtain a value of $A=2\cdot 10^{-12}$ m (Pa s)$^{-1}$. This value is in line with the values reported for  commercial membranes, which are in the range $10^{-14}$ to $10^{-10}$ (Pa s)$^{-1}$ \citep{Ruiz-Garcia2019,Lee1981,Dravzevic2014}.

{
\renewcommand\arraystretch{2}
\begin{table}[htpb]
\begin{center}
\resizebox{0.6\textwidth}{!}{%
\begin{tabular}{ c c  c }
\toprule
Symbol &  definition/value & units \\
 \midrule
	$k$ & $ 10^{-16}$ & m$^2$  \\
	$\epsilon$ & $ 0.7$ & - \\ 
    $ \phi_{a,0} $ & $ 35$ & g/m$^3 $ \\ 
    $ \phi_{b,0} $ & $ \phi_{a,0}$ &  g/m$^3$\\
    $ u_{in,av} $ & $ 0.1$ & m/s  \\
    $ D $ & $ 0.003$ & m  \\
    $ \ell $ & $ 0.0001$ & m  \\    
    $V_{s} $ & $27\cdot10^{-6}$ & m$^3$/mol \\
    $u_{av}$ & $0.1$ & m/s \\
    $\kin $ & $ \{10^{-10},\,\,10^{-5},\,\,10^{-2},\,\,10^{-1}\}$ & m$^3$/mol\\
    $\rho$ & 1000 & kg/m$^3$  \\
    $p_{out}$ & 1000 & kPa \\
    $p_{perm}$ & 0 & kPa  \\
    $\mu$ & $10^{-3}$ & Pa s \\
\bottomrule     
\end{tabular}
}
\caption{Summary of the numerical inputs for the physical quantities used in the simulations. Note that the units of the kinetic constant $\kin$ depend on the fact that we considered a binary reaction, whereas for the permeability $k$ and porosity $\epsilon$ we are considering the initial value (i.e., the value at $t=0$ h. }
\label{tab:consts}
\end{center}
\end{table}
}

Our model is able to include the scaling of the membrane given by the chemical reaction, which can modify the membrane permeability through the precipitation of a solid phase obtained as a product of the reaction. In this work we considered a range of kinetic reactions, going from very slow to fast reactions, i.e. with a value of the kinetic constant spanning four orders of magnitude, from $10^{-10}$ to $10^{-1}$ m$^3$/mol. We considered such a large range of kinetics since our main goal is not to focus on a specific system (and reaction) but to give a general description that can be applied to different specific situations.  

\section{Results}

\begin{figure*}
 \begin{subfigure}[b]{0.5\textwidth}
    \centering
   \includegraphics[width=\textwidth]{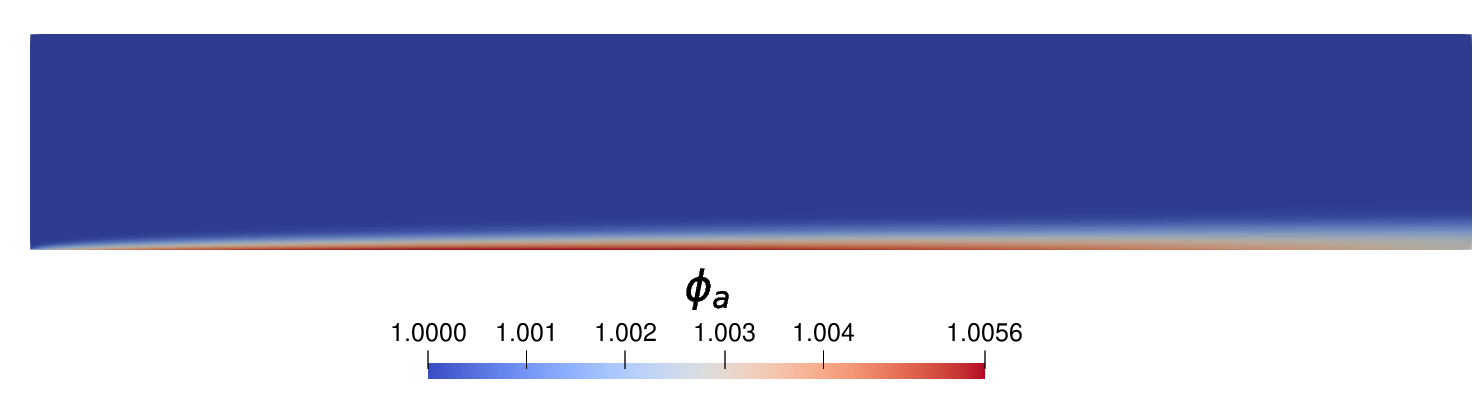}
   \caption{$\kin = 10^{-10}$ m$^3$/mol}
   \label{fig:lowk}
 \end{subfigure}
  \begin{subfigure}[b]{0.5\textwidth}
    \centering
   \includegraphics[width=\textwidth]{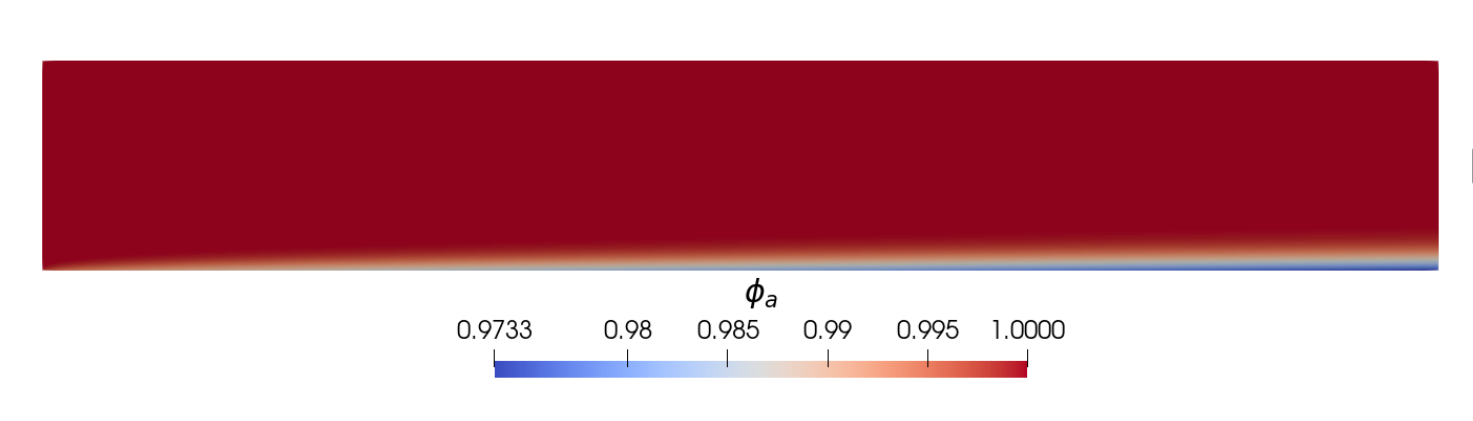}    
   \caption{$\kin = 10^{-1}$ m$^3$/mol}
   \label{fig:highk}
 \end{subfigure}
 \caption{Contour plot of the concentration within the channel given in units of the initial concentration $\tilde{\phi}_a=\phi_a \phi_{in,0} $. On the left: results reported for the lowest kinetic constant. On the right, results are reported for the highest kinetic constant. Note that the starting point of the legend is not zero and is different between the pictures to make the results more clear. }
    \label{fig:cplots}
\end{figure*}

In this work, we employ a fixed flow profile at the inlet, which when considering the properties in \cref{tab:consts}, gives a Reynolds number equal to 300. Therefore, the system operates in a fully developed laminar flow regime.
%
%
The first property that can be derived from this model is the polarisation of the membrane, which represents the accumulation of the solute at the interface of the membrane on the feed side. This is an undesired effect since it increases the osmotic pressure reducing the extraction of the permeate per unit of energy consumed in the process. 
We report the variation of the concentration profile in the domain in \cref{fig:cplots} for the lowest and highest chemical kinetic rate considered.
We can observe in \cref{fig:cplots} that the concentration at the membrane is different from the one in the bulk region in both cases. However, while the case with $\kin = 10^{-15}$ m$^3$/mol shows a higher concentration with respect to the bulk (see \cref{fig:lowk}), the case with the highest value of the kinetic rate ($\kin = 10^{-1}$ m$^3$/mol, see \cref{fig:highk}) shows a concentration smaller than the one in the bulk.

The latter observations show that the possible behaviour of the solution near the membrane strongly depends on the reaction kinetics. In the first case (the one represented in \cref{fig:lowk} corresponding to the lowest kinetic rate considered $\kin = 10^{-15}$ m$^3$/mol), we can observe the ``standard'' effect of the \textit{polarisation} of the membrane. During the filtration process, there is an accumulation of the solutes molecule on the feed side of the membrane, which results in a higher concentration of the solution at the membrane itself. On the opposite side, when the reaction rate is almost negligible (as for the case of $\kin = 10^{-15}$ m$^3$/mol), the solutes are now consumed by the reaction and they accumulate at the membrane interface, leading to the concentration profile observed in  \cref{fig:lowk}. In the latter case, the solute is now consumed almost instantly at the membrane interface. This result in a transport (convection and diffusion) limited profile of the concentration near the interface. 

The two opposite effects just described for the profile of the concentration at the membrane interface give an interesting effect on the evolution of the porosity and permeability profiles across the membrane. 
As the reaction proceeds, a new solid phase is formed which precipitates on the membrane modifying its structure and therefore its fluid dynamical behaviour. In particular, the solid phase generated during the reaction clogs the pores of the membrane, resulting in a variation of the porosity of the membrane with time. On top of this, since the concentration along the membrane (in the x-direction) decreases, we can expect a decrease in the overall reaction rate (which is proportional to the concentration) and therefore a difference in the permeability and porosity over the membrane. 
When we instead observe polarisation (i.e., in the case of the lowest reaction rate) the concentration near the membrane increases with the distance along the membrane direction. 
 
Therefore, we can expect that the reaction velocity (which depend on $\kin$ and the concentration), at a fixed time, will increase along the membrane interface for the case with polarisation (low kinetic reaction) and decrease along the membrane for high value of the kinetic reaction rate. We will show quantitatively these effects in the next section.

The description of the properties of the membrane (i.e., porosity, permeability, velocity through the membrane) can be given in terms of global quantities, that is to say quantities averaged over all the membrane length, which therefore becomes a function of time only. We will start our analysis by giving an account of these global properties. 

\begin{figure}
    \centering
   \includegraphics[width=.8\textwidth]{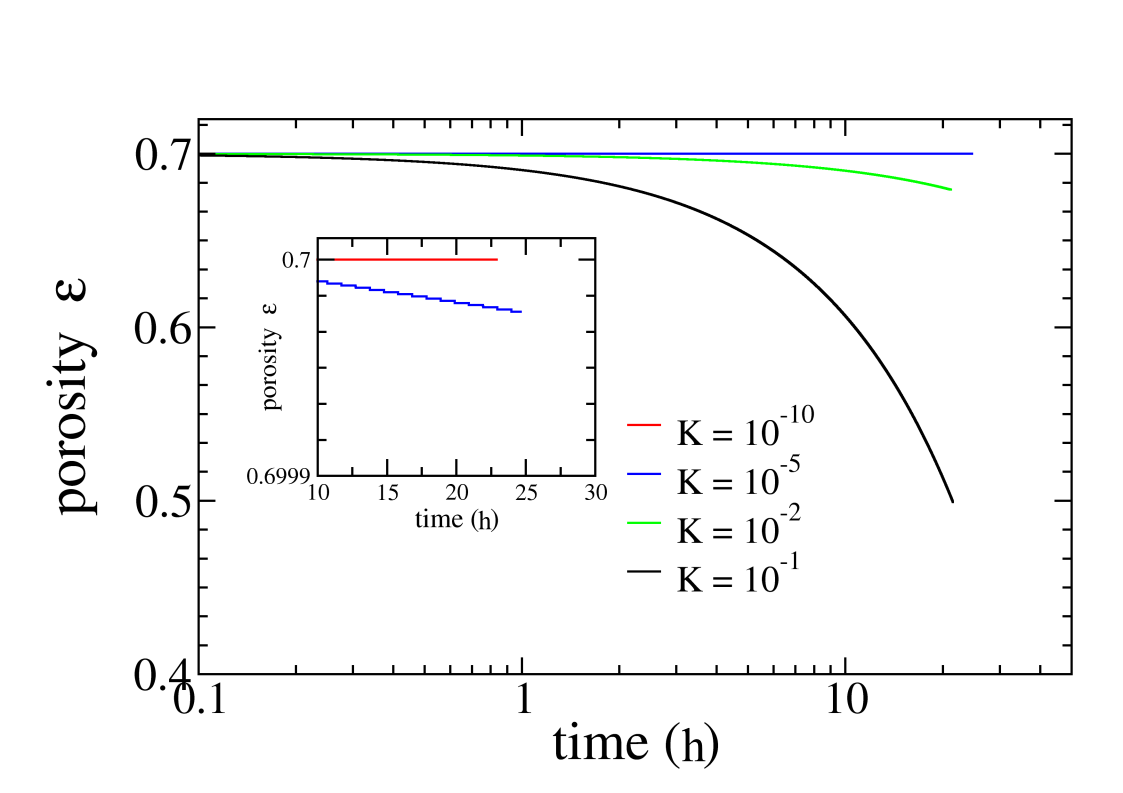}    
 \caption{Plot of the porosity versus time for all the systems considered, $\kin$ is in mol/m$^3$.   }
    \label{fig:eps}
\end{figure}

In \cref{fig:eps} we report the variation of the average of the porosity across the membrane as function of time. 
For the lowest kinetic reaction time considered there is no appreciable variation of the porosity after more than one day of operations. By increasing the kinetic reaction rate we can start to observe some deviations. In particular, for the highest value of the reaction rate the porosity decays to 60\% of its original value after only one day of operation.
This latter kind of results can be useful in determining the operation time that we can expect from a membrane given a certain composition of the feed.

\begin{figure}
    \centering
   \includegraphics[width=.8\textwidth]{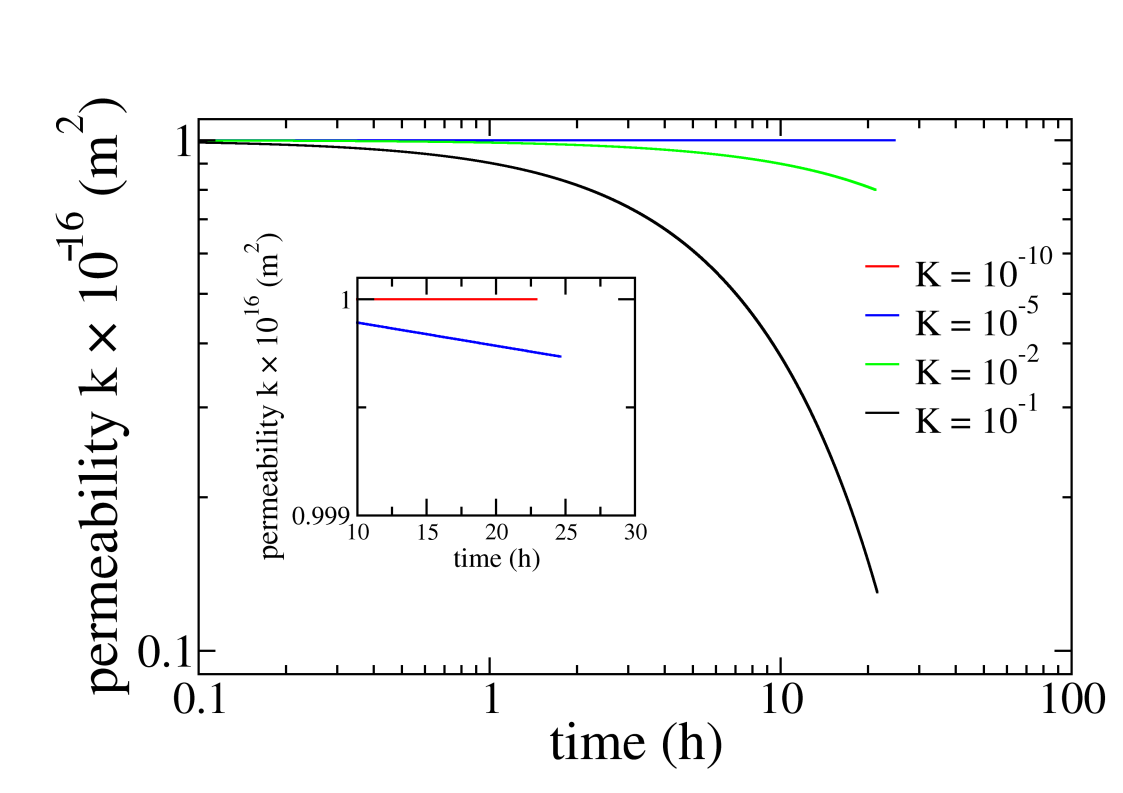}    
 \caption{Plot of the permeability as a function of time for all the systems considered $\kin$ is in mol/m$^3$.   }
    \label{fig:perm}
\end{figure}

The law of variation of the permeability with the reaction depends on the variation of the porosity, and in fact we can expect a similar behaviour. We reported the results for $k$ in \cref{fig:perm}, where we can see that there is an order of magnitude difference between the initial value of the permeability at time $t=0$ and after 28 h of operations.

The average velocity through the membrane obtained with the conditions specified is 1.8 $\mu$s, which decreases with the decrease of the permeability of the membrane up to 0.13 $\mu$s for the lowest value of $k$ and $\epsilon$ shown in \cref{fig:perm,fig:eps}.
In order to maintain the flow across the membrane in the given conditions of cross-flow in the channel and for the permeability and porosity given, we have to apply a pressure of approx 1800 kPa, which is needed to overcome an osmotic pressure of 1000 kPa, which reduces to 978 kPa in the case of the highest reaction rate. The difference in the osmotic pressure for the case $\kin=10^{-1}$ depends on the fact that for this case the concentration at the membrane is lower than the bulk (and lower than the case with polarisation) because of the very fast reaction (see \cref{fig:highk}). 

Despite the fact the osmotic pressure is smaller for the fastest reaction case, this case remains the worse in terms of permeate extraction, because the fast scaling of the membrane reduces quickly the porosity until the flux stops completely, i.e., we reach a value of $v= 0.13$ $\mu$s when we consider a fast reaction rate, against a value one order of magnitude higher for the case of the lowest reaction rate where we do not observe the scaling in the simulated time.

\subsection{Local profiles}

\begin{figure*}
 \begin{subfigure}[b]{0.5\textwidth}
    \centering
   \includegraphics[width=\textwidth]{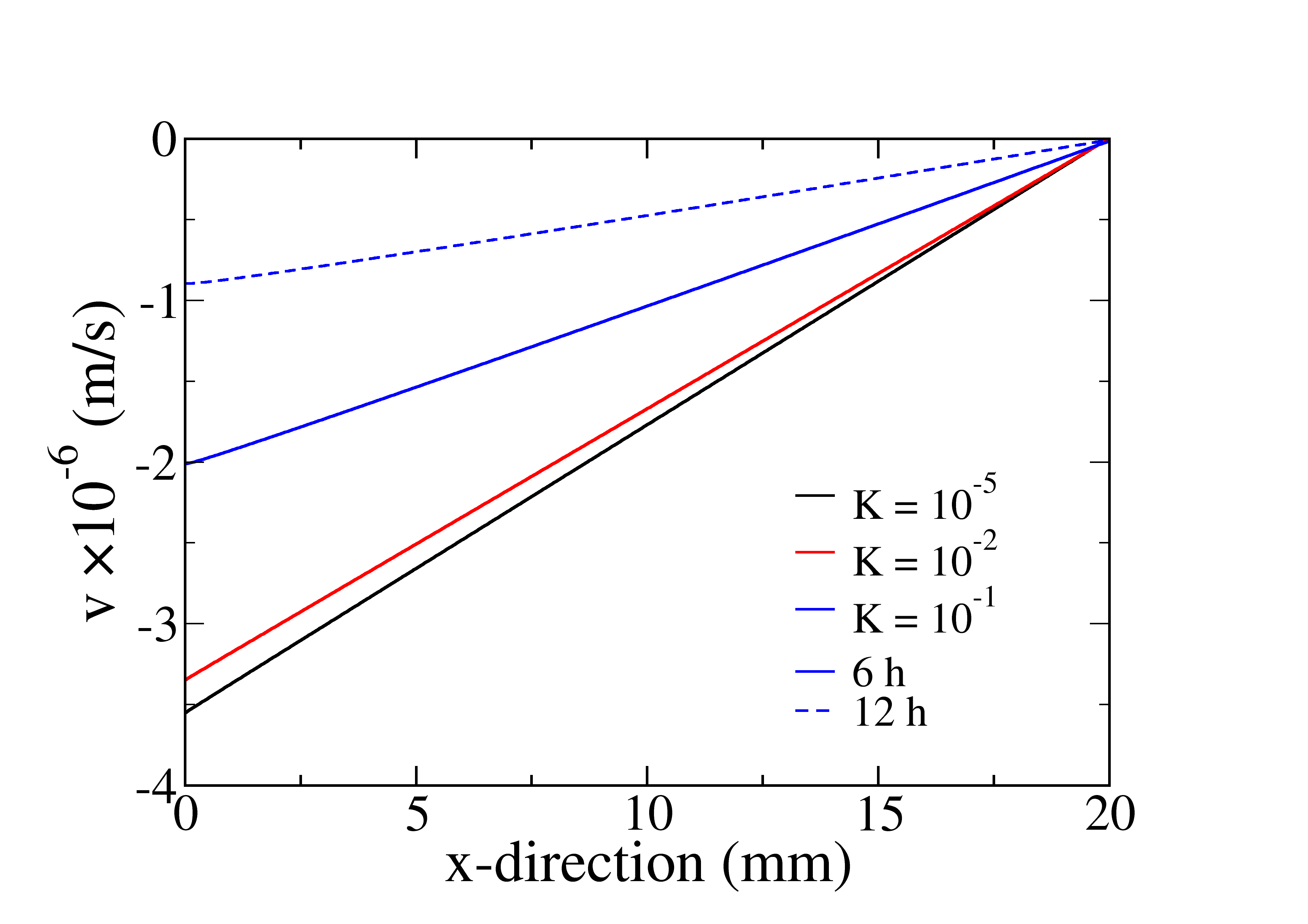}
    \caption{ $v$} 
     \label{fig:v-ax} 
    \vspace{2ex}
  \end{subfigure} 
    \begin{subfigure}[b]{0.5\textwidth}
   \centering
     \includegraphics[width=0.9\textwidth]{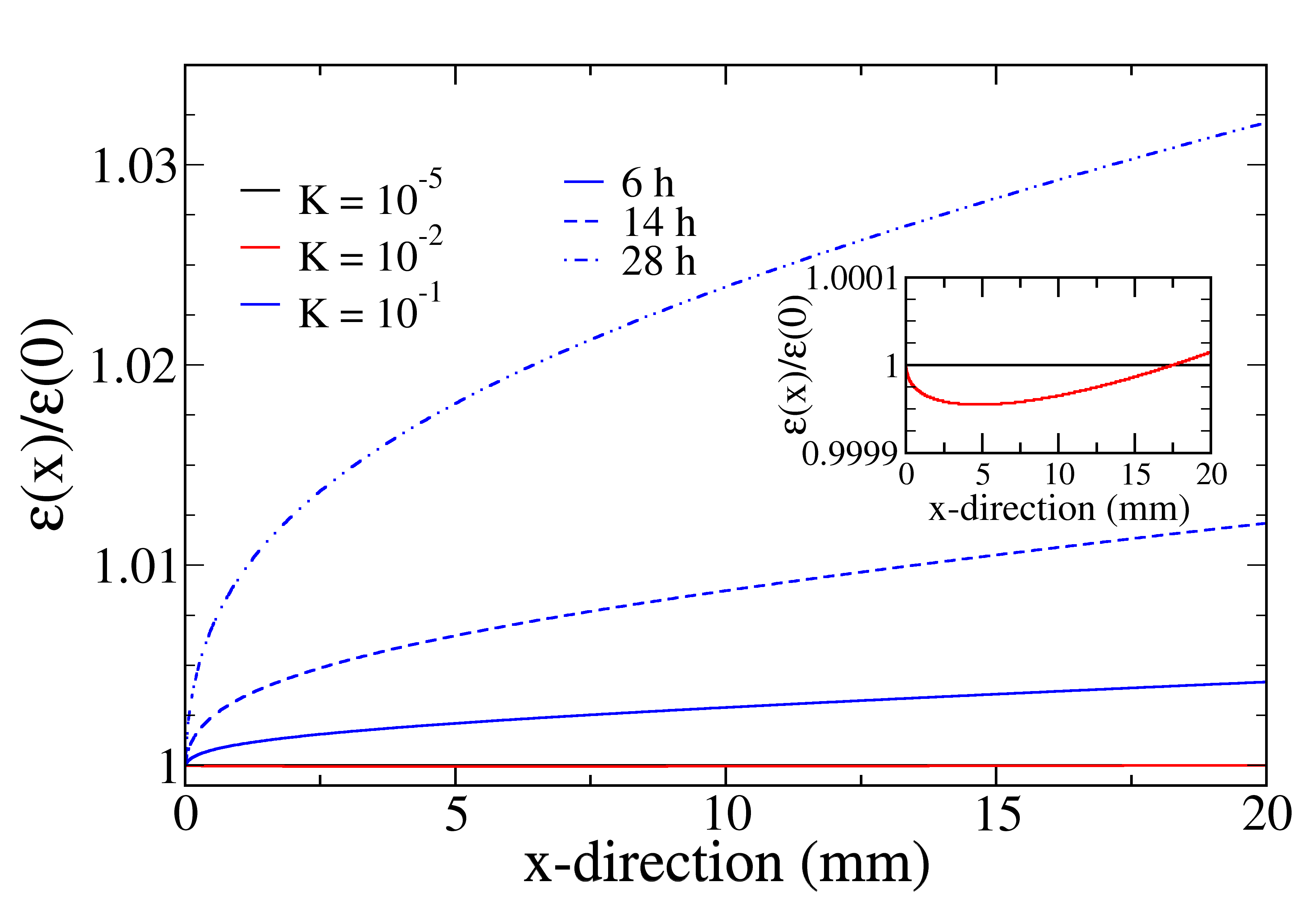}
    \caption{ porosity} 
    \label{fig:eps-ax} 
    \vspace{2ex}
  \end{subfigure}
 \caption{ Velocity through the membrane and porosity profile along the membrane at different times for different value of the kinetic reaction values. The continuous lines are the results at 6 h, the dashed lines at 12 h and the dot-dashed lines at 28 h. $\kin$ is in mol/m$^3$. }
    \label{fig:axial}
\end{figure*}

Local profiles at the membrane are analysed for the same range of parameters.  Since the smaller value of the kinetic constant ($\kin=10^-15,10^-10$ m$^3$/mol) gives the same behaviour in the time scale considered, we are showing only results for $\kin \geq 10^-5$ m$^3$/mol. We summarise our finding in \cref{fig:axial} where we reported the component $v$ of the velocity of the fluid (i.e., the velocity through the membrane) and the porosity are reported.

According to our analysis in the preceding sections, the flux along the membrane depends on different contributions. The first is the frictional pressure drop along the channel, which can cause considerable differences in the transmembrane pressure. Secondly, we have the polarisation effects, which increase the osmotic pressure (as it is proportional to the concentration difference across the two sides of the membrane), and finally, the scaling, which changes the permeability of the membrane itself. While the first contribution can be mitigated by a better design of the membrane modules \citep{Krawczyk2014}, the contributions of the last two effects are difficult to quantify \textit{a priori}, as it can be seen from \cref{fig:v-ax}. While the frictional pressure drop along the channel acts in all the systems in the same way (we are considering the same geometry and the same initial conditions for the flow), that is not true for the polarisation and scaling effects. In particular, the systems with a lower concentration suffer from polarisation at the membrane, as shown in the previous section, which increases the osmotic pressure. The system with higher reaction kinetics instead, does not suffer from the polarisation of the membrane (and in fact, the osmotic pressure is lower than the case at lower reaction rates, see the previous section). However, the scaling of the membrane, combined with the pressure drop in the channel is now dominating, resulting in an overall smaller flux (see blue curves against red and black curves in \cref{fig:v-ax}). 

In the variation of the porosity profiles in the membrane, we can observe the qualitative analysis we discussed in the previous section. For the larger kinetic rates, the transport-limited boundary layer on the membrane is reducing the reaction rate along the membrane. This, in turn, gives a porosity that increases along the membrane, with a value of $\epsilon$ at the outlet of the membrane, after 28 hours, 3\% larger than the value at the inlet (see dot-dashed blue curve in \cref{fig:eps-ax}). For the lowest reaction rates instead, we obtain the opposite behaviour: the polarisation increases the overall velocity along the membrane, which result in a reduction of the porosity along the x-direction, even though the low velocity of the reactions results in a very small variation (see the red continuous curve in \cref{fig:eps-ax}). 
Near the outlet of the membrane, the porosity increases, as a result of the reduction of the polarisation at this point of the domain.

%
\FloatBarrier

\section{Conclusions}

In this work, we presented a comprehensive computational model to describe the solute dynamics ed evolution near a membrane for desalination processes. In particular, we included a model to treat the scaling of the membrane as solid precipitated following a (general) chemical reaction. We connected the accumulation of solids at the membrane with porosity and permeability as described by the Darcy theory of porous media. Following this approach, we were able to give a full explicit model to derive the dynamical evolution of the filtration process by specifying a few initial parameters (e.g., the property of the solution and the kinetics of the reaction). 

The membrane is described as a dynamic boundary condition for the fluid mechanics and solute transport equations, which are coupled together through the osmotic pressure term, and therefore the flow through the membrane. We implemented our model in the widely used software package for CFD calculations \openfoam, and performed simulations for a selected range of operating conditions. Results show how this model can be used to predict the decay in the flux through the membrane due to the accumulation of the precipitated solid originating from the chemical reaction.

The formulation presented here has two main advantages which make it flexible and powerful in treating polarization. First, the proposed formulation can address all the interconnections between the different mechanisms (fluid dynamics, solute evolution, chemical reaction, scaling, and fouling) which affect the membrane performance. The second advantage is that the model can be easily extended to include more complex geometries, or models for the osmotic pressure (such as the Pitzer model \citep{Pitzer1973,Khraisheh2019}), fluid flow conditions in the system, as well as  more complex reactions paths. 

\bibliographystyle{apsrev}
\bibliography{bibliography_20Feb20}

\end{document}